\documentstyle[12pt, amscd, amssymb]{amsart}

\newcommand{\cL}{{\cal L}}
\newtheorem{thm}{Theorem}[section]
\newtheorem{prop}[thm]{Proposition}
\newtheorem{lem}[thm]{Lemma}
\newtheorem{cor}[thm]{Corollary}
\newtheorem{defn}[thm]{Definition}  

\newtheorem{rem}[thm]{Remark}
 
\newcommand{\ZZ}{{Z\!\!\!Z}}

\newcommand{\NN}{{I\!\!N}}

\newcommand{\ess}{\trianglelefteq}
\newcommand{\Rad}[1]{{\mathrm{Rad}}\:(#1)}
\newcommand{\End}[1]{{\mathrm{End}}\:(#1)}
\newcommand{\Hom}[2]{{\mathrm{Hom}}\:(#1,#2)}
\newcommand{\Jac}[1]{{\mathrm{Jac}}\:(#1)}
\newcommand{\Ker}[1]{{\mathrm{Ker}}\:(#1)}

\newenvironment{proof}{\par\noindent{\bf Proof.}}{$\square$\par\bigskip}   

\begin{document}

\title[ ]{On Semilocal Modules and Rings}

\author[ ]{Christian Lomp \\[+1mm]
\mbox{
 \rm Heinrich Heine Universit\"at}\\
 D-40225 D\"usseldorf, Germany \\
 Lomp\@math.uni-duesseldorf.de}

\thanks{ I would like to thank Nyguen Viet Dung for bringing the R. Camps
and W. Dicks paper to my attention. Moreover I want to express my thanks to
Patrick F. Smith and Robert Wisbauer for their interest and helpful suggestions. }

\maketitle

\begin{abstract} 

It is well-known that a ring $R$ is semiperfect if and only if $_RR$ (or $R_R$) is 
a supplemented  module. Considering {\em weak supplements} instead 
 of supplements we show that {\em weakly supplemented} modules $M$ are
{\em semilocal} (i.e., $M/Rad(M)$ is semisimple) and that $R$ is a
semilocal ring if and only if $_RR$ (or $R_R$) is weakly supplemented.  In
this context the notion of {\em finite hollow dimension} (or {\em finite
dual Goldie dimension}) of modules is of interest and yields a natural
interpretation of the Camps-Dicks characterization of semilocal rings.
Finitely generated modules are weakly supplemented if and only if they
have finite hollow dimension (or are semilocal).
\end{abstract}

\begin{center} \section{Preliminaries}\end{center}

Let $R$ be an associative ring with unit and throughout the paper $M$ will be a left unital
$R$-module.  By $N \ess M$ we denote an essential submodule $N\subset M$. 
 $M$ is {\it uniform} if $M\neq 0$ and every non-zero
submodule is essential in $M$, and  $M$ has {\it finite uniform
dimension} (or {\it finite Goldie dimension})  if there exists a sequence
\[\begin{CD}
0 @>>> \bigoplus_{i=1}^n U_i @>f>> M ,
\end{CD}\]
where all the $U_i$ are uniform and the image of $f$ is essential in $M$.
Then $n$ is called the uniform dimension of $M$ and we write $udim(M)=n$. 
 It is well known that this is equivalent to
 $M$ having no infinite independent family of non-zero submodules 
(there is a maximal finite independent family of uniform submodules). 

We denote a small submodule $N$ of $M$ by $N \ll M$. 
A module $M$ is said to be {\it hollow} if
$M\neq0$ and every proper submodule is small in $M$.  $M$ is said to have
{\it finite hollow dimension} (or {\it finite dual Goldie dimension}) if
there exists an exact sequence
\[\begin{CD}
M @>g>> \bigoplus_{i=1}^n H_i @>>> 0
\end{CD}\]
where all the $H_i$ are hollow and the kernel of $g$ is small in $M$. Then
$n$ is called the hollow dimension of $M$ and we write $hdim(M)=n$. 

B. Sarath and K. Varadarajan showed in \cite[Theorem 1.8]{sarath}
that in this case $M$ does not allow an epimorphism to a direct sum with
more than $n$ summands.  Dual to the notion of an independent family of
submodules we have:\\

\begin{defn}
	Let $M$ be an $R$-module and $\{ K_\lambda \}_\Lambda$ a family of
	proper submodules of $M$. $\{ K_\lambda \}_\Lambda$ is called 
	{\it coindependent} (see \cite{takeuchi}) 
	if for every $\lambda \in \Lambda$ and finite subset 
	$J \subseteq \Lambda \setminus
	\{\lambda \}$ $$ K_\lambda + \bigcap_{j \in J} K_j = M$$
	holds(convention: if $J$ is the empty set, then set $\bigcap_J K_j
:= M$). 
\end{defn}

 For finitely generated modules it usually suffices to consider coindependent
families of finitely generated submodules as our next observation shows.

\begin{lem} \label{lemma_fg}
Let $M$ be a finitely generated $ R$-module and $\{ N_1, \ldots , N_m \}$
a coindependent family of submodules. Then there exists a coindependent family of
 finitely generated submodules $L_i \subseteq N_i$,  $1\leq i \leq m$. 
\end{lem}

\begin{proof}
Since $M$ is finitely generated, for each $1 \leq i \leq m$, there exist
finitely generated submodules $X_i \subseteq N_i$ and $Y_i \subseteq
\bigcap_{j\neq i} N_j$ such that $ X_i + Y_i=M$.
Let $ L_i := X_i + \sum_{j\neq i} Y_j \subseteq N_i$. 
 As $L_i + \bigcap_{j\neq i} L_j \supseteq X_i + Y_i = M$ holds
the result follows.
 \end{proof}

\begin{thm}[\small{Grezeszcuk, Puczy{\l}owski, Reiter, Takeuchi,
Varadarajan}] 
 For an $R$-module $M$ the following statements are equivalent:
 \begin{enumerate}
 	\item[(a)] $M$ has finite hollow dimension.
	\item[(b)] $M$ does not contain an infinite coindpendent family 
		of submodules.
	\item[(c)] There exists a unique number $n$ and a coindependent 
		family $\{ K_1, \ldots , K_n \}$ of proper submodules, 
		such that $M/K_1, \ldots, M/K_n$ are hollow modules
		and $K_1 \cap \cdots \cap K_n \ll M$.
	\item[(d)] For every descending chain 
		$K_1 \supset K_2 \supset K_3 \supset \cdots$ 
		of submodules of $M$, there exits 
		 a number $n$ such that $K_n/K_m \ll M/K_m$, for all $m \geq n$.
 \end{enumerate}
\end{thm}

\begin{proof}
The equivalence of $(b),(c),(d)$ can be found in \cite{gr}. The
equivalence of $(a)$ and $(c)$ is given by the chinese remainder theorem
(see \cite[9.12]{wisbauer}). 
\end{proof}

\begin{rem}\label{hdim_prop}\em
Let $M$ be an $R$-module and $N, L$ submodules of $M$. Then the following
properties hold:
\begin{enumerate}
	\item[(1)] $hdim(M/N) \leq hdim(M)$;
	\item[(2)] $N \ll M \Rightarrow hdim(M) = hdim(M/N)$;
	\item[(3)] $hdim(N\oplus L) =\ hdim(N) + hdim(L)$.
\end{enumerate}
Moreover if $M$ is self-projective and has finite hollow dimension, then
every surjective endomorphism is an isomorphism.

We refer to \cite{gr}, \cite{hanna}, \cite{herbera},\cite{lomp},
\cite{reiter}, \cite{takeuchi} and \cite{var1} for more information on
dual Goldie dimension.
\end{rem}

The following theorem can be seen as an attempt to transfer R. Camps and
W. Dicks characterization of semilocal rings \cite{camps} to
arbitrary modules with finite hollow dimension. Denote by
$\cL (M)$ the lattice of submodules of a module $M$. 

\begin{thm} \label{mod_Camps_Dicks}
  For $M$ the following statements are equivalent:
 \begin{enumerate}
  \item[(a)] $M$ has finite hollow dimension.
  \item[(b)] There exists an $n\in \NN$ and a mapping
             $d:\cL (M) \rightarrow \{ 0,1, \ldots , n \}$ such that
             for all $N,L \in \cL (M)$:
             \begin{enumerate}
               \item[(i)] If $d(N)=0$, then  $N=M$.
               \item[(ii)] If $N+L=M$, then $d(N \cap L) = d(N) + d(L)$.
             \end{enumerate}
  \item[(c)] There exists a partial ordering $(\cL (M) , \leq)$ such
	     that 
	     \begin{enumerate}
		\item[(i)]  $(\cL (M), \leq)$ is an artinian poset;
		\item[(ii)] for all $N,L \in \cL (M)$ with $N+L=M$:
			    if $L\neq M$, then $N > N \cap L$.
	     \end{enumerate}
 \end{enumerate}
\end{thm}

\begin{proof}
$(a) \Rightarrow (b)$ Let $d(N):=hdim(M/N)$;
     then the conditions $(i)$ and $(ii)$ are easily checked.

$(b) \Rightarrow (c)$ Let $N < L :\Leftrightarrow d(N) < d(L)$ and $N=L
     :\Leftrightarrow d(N) = d(L)$
     then $(\cL (M) ,\leq)$ is artinian.
     Let $N, L \in \cL (M)$ with $N+L=M$ and $L\neq M$. By $(i)$ and
   $(ii)$ we have $d(N \cap L)= d(N) + d(L) < d(N)$. Hence $N > N\cap L$.

$(c) \Rightarrow (a)$ Assume that $\{ K_i \}_\NN$ is an infinite
     coindependent family of submodules of $M$. Then we have for all
     $i\in \NN$: $K_1 \cap \cdots \cap K_i + K_{i+1} = M$ and $K_{i+1}\neq M$.
     Hence by $(ii)$ we get the infinite descending chain
     $$ K_1 > K_1 \cap K_2 > \cdots > K_1\cap \cdots \cap K_i > \cdots $$
     contradicting property $(i)$. Hence $M$ does not
     contain an infinite coindependent family of submodules.
\end{proof}

\begin{center} \section{Weakly Supplemented Modules}\end{center}

Dual to a complement of a submodule $N$ of $M$ the {\it supplement} of $N$
is defined as a submodule $L$ of $M$ minimal with respect to $N+L=M$. This
is equivalent to $N+L=M$ and $N\cap L \ll L$. Recall that $M$ is {\it
supplemented } if every submodule has a supplement in $M$. 

 More generally, a submodule
$N$ of $M$ has a {\it weak supplement} $L$ in $M$ if $N+L=M$ and $N \cap L\ll M$,  
 and $M$ is called {\it weakly supplemented} if every
submodule $N$ of $M$ has a weak supplement (see Z\"oschinger \cite{zoeschinger78a}).
 Examples for weakly supplemented modules are supplemented, artinian, linearly compact,
uniserial or hollow modules.
 For supplemented modules over commutative local noetherian
rings we refer to \cite{zoeschinger78a}, \cite{zoeschinger78b},
\cite{zoeschinger86} and \cite{rudlof}. 

 Before we give a summarizing list of properties of weakly supplemented
 modules, we will state a general result:\\

 \begin{prop} \label{gen_ws} 
  For a proper submodule $N\subset M$, the following are equivalent: 

  \begin{enumerate}
	\item[(a)] $M/N$ is semisimple;

	\item[(b)] for every $L \subseteq  M$ there exists a submodule
		$K\subseteq M$ such that $L+K=M$ and $L\cap K \subseteq N$;

	\item[(c)] there exists a decomposition $M=M_1 \oplus M_2$ such that
		$M_1$ is semisimple, $N \ess M_2$ and $M_2/N$ is semisimple.
  \end{enumerate}
\end{prop}

\begin{proof}
(a) $\Rightarrow $(c)
   Let $M_1$ be a complement of $N$.  
   $M_1 \simeq (M_1 \oplus N)/N$ is a direct summand in $M/N$, 
   hence semisimple and there is a semisimple submodule $M_2/N$
   such that $(M_1 \oplus N)/N \oplus M_2/N=M/N$. Thus $M=M_1+M_2$ and
   $M_1\cap M_2 \subseteq N \cap M_1 = 0$ implies $M=M_1 \oplus M_2$.
   Since $M_1$ is a complement we  have by the natural isomorphisms 
   $N \simeq (M_1 \oplus N)/M_1 \ess M/M_1\simeq M_2$ that $N \ess M_2$.

(c) $\Rightarrow $(a) $\Rightarrow$ (b) clear.

(b) $\Rightarrow $(a) Let $L/N \subseteq M/N$;
   then there exists a submodule $K\subseteq M$
   such that $L+K=M$ and $L\cap K \subseteq N$.
   Thus $L/N \oplus (K+N)/N = M/N$. Hence every submodule of $M/N$ is a
   direct summand.
\end{proof}

Let $\Rad{M}$ denote the radical of $M$.  We call $M$ 
a {\it semilocal module} if $M/\Rad{M}$ is semisimple. Any semilocal module $M$ 
is a {\it good} module, i.e.,  for every 
homomorphism $f:M\rightarrow N$, $f(\Rad{M}) = \Rad{f(M)}$  (see \cite{wisbauer}). 

We call $N$ a {\it small cover} of a module $M$ if there
exists an epimorphism $f:N \rightarrow M$ such that $\Ker{f} \ll M$.
 Then $f$ is called a {\it small epimorphism}.  $N$ is called a {\it
flat cover}, {\it projective cover} resp. {\it free cover} of $M$ if $N$
is a small cover of $M$ and $N$ is a flat, projective resp. free module.
Note that this definition of a flat cover is different from Enochs'
definition.

\begin{prop} \label{ws}
Assume $M$ to be weakly supplemented. Then: 
 \begin{enumerate}
    	  \item[(1)] $M$ is semilocal;
    	  \item[(2)] $M=M_1 \oplus M_2$ with $M_1$ semisimple, $M_2$
		     semilocal and $\Rad{M} \ess M_2$;
	  \item[(3)] every factor module of $M$ is weakly supplemented;
          \item[(4)] any small cover of $M$ is weakly supplemented;
	  \item[(5)] every supplement in $M$ and every direct summand 
		     of $M$ is weakly supplemented. 
  \end{enumerate}
\end{prop}

\begin{proof}
(1) and (2) follow from Proposition \ref{gen_ws} since for every 
	$L\subseteq M$ there exists a weak supplement $K\subseteq M$ 
	such that $L+K=M$ and $L\cap K \subseteq \Rad{M}$.

(3) Let $f : M \rightarrow N$ be an epimorphism and $K\subset N$, then
 $f^{-1} (K)$ has a weak supplement $L$ in $M$ and it is straightforward 
 to prove that $f (L)$ is a weak supplement of $K$ in $N$.

(4) Let $N$ be a small cover of $M$ and $f: N \rightarrow M$ be a small epimorphism.
First note that $f^{-1} (K) \ll N$ for every $K \ll M$ holds since  
$\Ker{f} \ll N$. Let $L \subset N$. Then $f (L)$ has a weak
supplement $X$ in $M$. Again it is easy to check that $f^{-1}(X)$
is a weak supplement of $L$ in $N$.

(5) If $N \subseteq M$ is a supplement of $M$, then $N + K = M$ for some
   $K \subseteq M$ and $K\cap N \ll N$. By (3),
   $M/K \simeq N/(N\cap K)$ is weakly supplemented and by (4), $N$ is
   weakly supplemented. Direct summands are supplements and
   hence weakly supplemented.
\end{proof}

Let $length(M)$ denote the length of the module $M$. 

\begin{cor} \label{semiprimitive_ws}
  An $R$-module $M$ with $\Rad{M}=0$ is weakly supplemented if and only
  if $M$ is semisimple. In this case $hdim(M)=length(M)$ holds.
\end{cor}

\begin{proof}
This follows by Proposition \ref{ws}{(1)}.
\end{proof}       
 
We need the following technical lemma to show that every finite sum
of weakly supplemented modules is weakly supplemented.

\begin{lem}
        Let $M$ be an $R$-module with submodules $K$ and $M_1$. Assume
	$M_1$ is weakly supplemented and $M_1 + K$ has a weak supplement
	in $M$. Then $K$ has a weak supplement in $M$.
\end{lem}

\begin{proof}
   By assumption $M_1 + K$ has a weak supplement $N \subseteq M$, such
   that $M_1 + K + N =M$ and $(M_1 + K) \cap N \ll M$. Because $M_1$ is
   weakly supplemented, $(K + N)\cap M_1$ has a weak supplement
   $L\subseteq M_1$. So
   $$ M=M_1 + K + N = L + ((K+N)\cap M_1) + K + N = K + (L +N) \mbox{ and } $$
   $$K \cap (L+N) \subseteq ((K+L)\cap N) + ((K+N) \cap L
   \subseteq ((K+M_1)\cap N) + ((K+N)\cap L) \ll M.$$
   Hence $N + L$ is a weak supplement of $K$ in $M$. 
\end{proof}

\begin{prop}
 Let $M=M_1 + M_2$, where $M_1$ and $M_2$ are
	weakly supplemented, then $M$ is weakly supplemented.
\end{prop}

\begin{proof}
   For every submodule $N \subseteq M$, $M_1 + (M_2 + N)$ has the trivial
   weak supplement $0$ and by the Lemma above $M_2 + N$ has a weak 
   supplement in $M$ as well. 
   Applying the Lemma again we get a weak supplement for $N$.
\end{proof}

\begin{cor}\label{fsum_ws}
   Every finite sum of weakly supplemented modules is weakly supplemented.
\end{cor}

The relationship between the concepts 'hollow dimension' and
'weakly supplemented' is expressed in the following theorem.

\begin{thm} \label{fg_ws} Consider the  following properties:
   \begin{enumerate}
     \item[(i)] $M$ has finite hollow dimension;
     \item[(ii)] $M$ is weakly supplemented;
     \item[(iii)] $M$ is semilocal.
   \end{enumerate}
   Then $(i) \Rightarrow (ii) \Rightarrow (iii)$ and $hdim(M) \geq
   length(M/Rad(M))$ holds. \\
   If $Rad(M) \ll M$ then $(iii) \Rightarrow (ii)$ holds.\\
   If $M$ is finitely generated then $(iii) \Rightarrow (i)$ and
   $hdim(M)=length(M/Rad(M))$ holds.
\end{thm}

\begin{proof}
 $(i) \Rightarrow (ii)$ There is a small epimorphism
    $f : M \rightarrow \bigoplus_{i=1}^n H_i$ with hollow modules $H_i$.
    Since hollow modules are (weakly) supplemented we get by
    Corollary \ref{fsum_ws} 
    that $\bigoplus_{i=1}^n H_i$ is weakly supplemented. Since $f$ is a
    small epimorphism we get by Proposition \ref{ws}(4) that $M$ is weakly
    supplemented.

 $(ii) \Rightarrow (iii)$ by Propositon \ref{ws}(1).

 If $Rad(M) \ll M$, then $(iii) \Rightarrow (ii)$
    follows by Proposition \ref{ws}(4).

 If $M$ is finitely generated and $(iii)$ holds,
    then $M$ is a small cover of $M/\Rad{M}$.
    By Corollary \ref{semiprimitive_ws}, 
    $hdim(M/\Rad{M}) = length(M/\Rad{M})$, and by
    remark \ref{hdim_prop}(2), $hdim(M) = length(M/\Rad{M})$.  
\end{proof}

\begin{center} \section{Semilocal Modules and Rings} \end{center}

Let $Gen(M)$ denote the class of $M$-generated modules.
\smallskip

\begin{thm}\label{semilocal_modules}
 The following statements about $M$ are equivalent:
  \begin{enumerate}
	\item[(a)] $M$ is semilocal;
	\item[(b)] any $N\in Gen(M)$ is semilocal;
	\item[(c)] any $N\in Gen(M)$ is a direct sum of a semisimple
		   module and a semilocal module with essential radical;
	\item[(d)] any $N\in Gen(M)$ with small radical is weakly
		   supplemented;
	\item[(e)] any finitely generated $N \in Gen(M)$ has finite hollow
		   dimension.
 \end{enumerate}
\end{thm}

\begin{proof}
$(a) \Rightarrow (b)$ For every $N\in Gen(M)$ there exists a set
$\Lambda$ and an epimorphism $f: M^{(\Lambda)} \rightarrow N$. 
Since $f(\Rad{M^{(\Lambda)}}) \subseteq \Rad{N}$ and
$M^{(\Lambda)}/\Rad{M^{(\Lambda)}} \simeq (M/\Rad{M})^{(\Lambda)}$ always
holds we get an epimorphism ${\bar f}: (M/\Rad{M})^{(\Lambda)} \rightarrow
N/\Rad{N}$. Hence $N$ is semilocal.

$(b) \Rightarrow (a)$ trivial.

$(b) \Leftrightarrow (c)$ by Proposition \ref{gen_ws}.

$(b) \Leftrightarrow (d) \Leftrightarrow (e)$ by Theorem \ref{fg_ws}
\end{proof}

%

 
 Recall that the ring $R$ is {\it semilocal} if $_RR$ (or $R_R$) is a semilocal $R$-module. 

\begin{cor} \label{semilocal}
 For a ring $R$ the following statements are equivalent:
 \begin{enumerate}
  \item[(a)] $_RR$ is weakly supplemented;
  \item[(b)] $_RR$ has finite hollow dimension;
  \item[(c)] $R$ is semilocal;
  \item[(d)] $R_R$ has finite hollow dimension;
  \item[(e)] $R_R$ is weakly supplememted.
 \end{enumerate}
 In this case $hdim(_RR) = length(R/\Jac{R}) = hdim(R_R)$.
\end{cor}

\begin{proof}
Apply Theorem \ref{fg_ws} and use that 'semilocal' is a left-right
symmetric property.
\end{proof}

\begin{rem} \em
 Consider the ring 
      $$R:=\ZZ_{p,q} := \left \{ \frac{a}{b} | a,b \in \ZZ, b \neq 0,  p
\nmid b \mbox{ and } q \nmid b \right\},$$
      where $p$ and $q$ are primes. Then $R$ is a
      commutative uniform semilocal noetherian domain with two maximal ideals.
      Since $R$ is uniform, the decomposition of
      $R/\Jac{R}$ cannot be lifted to $R$. Moreover the maximal ideals $pR$ and $qR$ are
      weak supplements but not supplements of each other.
      So $R$ is a semilocal ring which is not semiperfect.         
\end{rem}

 For our next result we need the following: 

\begin{lem} \label{lem_sem1}
 Let $R$ be a ring, $r,a \in R$ and $b:=1-ra$. Then
 $Ra \cap Rb = Rab$.
\end{lem}

\begin{proof}
 $x\in Ra \cap Rb$, then $x = ta = sb = s(1-ra) \Rightarrow s = (t+sr)a
 \in Ra$. Hence $ Ra \cap Rb \subseteq Rab$. Conversely $Rab =
 Ra(1-ra) = R(1-ar)a \subseteq Ra \cap Rb$.
\end{proof}

We are now ready to give characterizations of semilocal rings in terms of
finite hollow dimension and to prove results from Camps-Dicks 
(see \cite[Theorem 1]{camps}) in a module-theoretic way. 

Note that for a semilocal ring $R$,  $_RR$ is a good
module, and so for any left $R$-module $N$ we have 
 $\Rad{M} = \Jac{R}\, M$ (see \cite[23.7]{wisbauer}). \\

\begin{thm}\label{thetheorem}
 For any ring $R$ the following statements are equivalent:
 \begin{enumerate}
  \item[(a)] $R$ is semilocal;
  \item[(b)] every left $R$-module is semilocal;
  \item[(c)] every left $R$-module is the direct sum of a semisimple
module and a semilocal module with essential radical;
  \item[(d)] every left $R$-module with small radical is weakly
supplemented;
  \item[(e)] every finitely generated left $R$-module has finite hollow
dimension;
  \item[(f)] every product of semisimple left $R$-modules is semisimple;
  \item[(g)] there exists an $n\in \NN$ and
        a map $d:R \rightarrow \{ 0,1, \ldots, n\}$ such that for
        all $a,b \in R$ the following holds:
        \begin{enumerate}
          \item[(i)] $d(a) = 0  \Rightarrow a$ is a unit;
          \item[(ii)] $d(a(1-ba)) = d(a) + d(1-ba)$;
        \end{enumerate}
  \item[(h)] there exists a partial ordering $(R,\leq)$ such that:
	\begin{enumerate}
	  \item[(i)] $(R, \leq)$ is an artinian poset;
	  \item[(ii)] for all $a,b \in R$ such that
        	$1-ba$ is not a unit, we have $a > a(1-ba)$.
	\end{enumerate}
 \end{enumerate}
 In this case $hdim(R) \leq n$ holds. 
\end{thm}

\begin{proof}
$(a) \Leftrightarrow (b) \Leftrightarrow  (c) \Leftrightarrow (d)
\Leftrightarrow (e)$ follow from Theorem \ref{semilocal_modules}.

$(b) \Rightarrow (f)$ By the remark above, semilocal rings are good rings
and hence $\Rad{M} = \Jac{R}M$ holds for every left $R$-module. Let $M$ be
a product of semisimple modules. Since for all $m\in M$ $\Jac{R}Rm=0$
holds as $Rm$ is semisimple, we have $\Rad{M}=\Jac{R}M = 0$.  By $(b)$ $M$
is semisimple. 

$(f) \Rightarrow (a)$  $R/\Jac{R}$ is a submodule of a
product of simple modules. By $(f)$ this product is semisimple and so is
$R/\Jac{R}$. 

$(a) \Rightarrow (g)$  By Corollary \ref{semilocal} $_RR$ has finite
hollow dimension. By Theorem \ref{mod_Camps_Dicks} there is a map $d'$.
Set $d(a):=d'(Ra)$ for all $a\in R$ and (i) and (ii) follow easily from
the properties of $d'$.

$(g) \Rightarrow (h)$
  Let $a < b :\Leftrightarrow d(a) < d(b)$ and 
 $a=b :\Leftrightarrow d(a)=d(b)$
  for all $a, b \in R$. If $1-ba$ is not a unit, then $d(1-ba) >0$
  implies $d(a) < d(a(1-ba))$ and hence $a > a(1-ba)$.

$(h) \Rightarrow (a)$
  Assume that there exists a left ideal $I\subset R$ that has no weak
  supplement.
  Then we can construct an infinite descending chain of elements
  $$ 1 > b_1 > b_2 > \cdots > b_n > \cdots $$
  such that for  all $n\in \NN$ we have $I + Rb_n = R$.
  Since $(R, \leq)$ is artinian - this is a contradiction, hence
  $I$ must have a weak supplement in $R$. By Corollary \ref{semilocal}
  $R$ is semilocal.\\
  We can construct the chain as follows:
  Let $n=1$. Since $I \not \ll R$ there is an $a \in I$ such that $1-a$ is
  not a unit in $R$. Hence $1 > 1-a =: b_1$ and $I + Rb_1 = R$ holds.\\
  Now assume that we constructed a chain
  $1 > b_1 > b_2 > \cdots > b_n$ for $n\geq 1$
  with $I + Rb_n = R$. By assumption
  $I \cap Rb_n \not \ll R$ implies that there is an $r\in R$,
  such that $rb_n \in I$ and $x:=1 - rb_n$ is not a unit in $R$. 
  Hence
  $$b_n > b_n (1 - rb_n) =b_n x =: b_{n+1} .$$
  Moreover, by the modularity law, we have
  $Rb_n = (I\cap Rb_n) + (Rb_n \cap  Rx)$.
  Together with Lemma \ref{lem_sem1}, 
  $R = I + Rb_n = I + (Rb_n \cap Rx) = I + Rb_{n+1}$ holds.
\end{proof}
\begin{rem}\em
	Theorem \ref{thetheorem} generalizes the well-known fact that a
	 ring $R$ is semiperfect if and only if every finitely generated
	$R$-module is supplemented. 
\end{rem}

Recall that every finitely generated $R$-module over a semiperfect ring 
$R$ has a projective cover.

\begin{cor}
	Every finitely generated $R$-module over a semilocal ring $R$ is a
	direct summand of a module having a finitely generated free cover.  
\end{cor}

\begin{proof}
Let $M$ be a finitely generated $R$-module. Then there exists a number $k$ and an epimorphism $f:R^k \rightarrow M$.
Since $R$ is semilocal, $R^k$ is weakly supplemented. Hence $K:=\Ker{f}$ has a weak supplement $L\subseteq R^k$.
Thus the natural projection $R^k \rightarrow M \oplus (R^k / L)$ with kernel $K\cap L \ll R^k$ implies that $R^k$ is a projective cover 
for $M \oplus (R^k/L)$.
\end{proof}

Comparing semiperfect and semilocal rings the following fact is of interest:

\begin{thm} For a ring $R$ the following statements are equivalent.
  \begin{enumerate}
    \item[(a)] $R$ is semiperfect;
    \item[(b)] $R$ is semilocal and every simple $R$-module has a flat cover;
    \item[(c)] $R$ is semilocal and every finitely generated $R$-module has a flat cover.
  \end{enumerate}
\end{thm}

\begin{proof} 
$(a) \Rightarrow (c)$ holds since projective modules are flat.
 
$(c) \Rightarrow (b)$ is trivial. 

$(b) \Rightarrow (a)$ Assume $R$ is semilocal and consider 
$R/\Jac{R} = E_1 \oplus \cdots \oplus E_n$ with $E_i$ simple $R$-modules. 
Every simple $R$-module is isomorphic to one of the $E_i's$. 
By hypothesis every $E_i$ has a flat cover $L_i$. 
Thus $L:=L_1 \oplus \cdots \oplus L_n$ is a flat cover of $R/\Jac{R}$. 
Hence we obtain the following diagram:
$$\begin{CD}
@.          R \\
@.        @VVV \\
L       @>f>> R/\Jac{R} @>>> 0 \\
\end{CD}$$
that can be extended by a homomorphism $g:R \rightarrow L$. 
Since $f$ is a small epimorphism and $gf$ is epimorph, 
$g$ must be epimorph with $\Ker{g} \subseteq \Ker{gf} = \Jac{R}$.
Hence $R$ is a projective cover of the flat module $L$. 
By \cite[36.4]{wisbauer}, $L \simeq R$ and hence all $L_i$ must be projective. 
Thus each simple $R$-module has a projective cover and so $R$ is semiperfect 
(see \cite[42.6]{wisbauer}).
\end{proof}
\begin{rem} \em
It follows also from Theorem \cite[36.4]{wisbauer} that a ring $R$ is
semisimple if and only if $R$ is semilocal and every simple $R$-module is
flat. Since in this case $R$ is a projective cover of the flat module
$R/\Jac{R}$ and hence $R \simeq R/\Jac{R}$ holds.  
\end{rem}

The following result was first proved by T.Takeuchi in \cite{takeuchi94}.
We will give a new proof of his result.

\begin{thm}[Takeuchi] \label{pr_hdim}
Let $M$ be a self-projective $R$-module. Then $M$ has finite hollow dimension 
if and only if $S:=\End{M}$ is semilocal. Moreover we have 
$hdim(_RM) = hdim(S)$. 
\end{thm}

\begin{proof}
$\Rightarrow: $ 
Let $\{ I_1, \ldots , I_n \}$ be a coindependent family of proper left ideals of
$_SS$. By Lemma \ref{lemma_fg}, we may assume that the ${I_k} 's$ are
finitely generated. Consider the epimorphism
\[\begin{CD} 
	S @>>> 	\bigoplus_{k=1}^n S/I_k @>>>0. 
\end{CD}\]
Applying $M\otimes_S -$  we get the exact sequence
\[\begin{CD}
	M @>>> \bigoplus_{k=1}^n M/MI_k @>>>0,
\end{CD}\]
since $M \otimes_S S/I_k \simeq M/MI_k$.
 We have $I_k = \Hom{M}{MI_k}$ and hence $MI_k \neq M$. 
 Thus $hdim(S) \leq hdim(_RM)$ and so $S$ is semilocal by Corollary \ref{semilocal}.

$\Leftarrow:$ Consider an epimorphism (with $N_i \neq M$)
\[\begin{CD} 
	M @>>> 	\bigoplus_{i=1}^n M/N_i @>>>0. 
\end{CD}\]
Since $M$ is self-projective,  $\Hom{M}{-}$ yields an exact  sequence
\[\begin{CD}
	S @>>> \bigoplus_{i=1}^n \Hom{M}{M/N_i} @>>> 0, 
\end{CD}\]
showing that $hdim(_RM) \leq hdim(S)$. 
\end{proof}

\begin{rem}\em 
 More generally, if $P$ is an $M$-projective module that generates $M$, then
one can apply $\Hom{P}{-}$ in the same way as in Theorem \ref{pr_hdim} to
obtain $hdim(_RM)  \leq hdim(_S\Hom{P}{M})$, where $S:=\End{P}$. 
 \end{rem}

The following Corollaries are immediate consequences from Takeuchi's
result.

\begin{cor}
A ring $R$ is semilocal if and only if every finitely generated,
self-projective left (or right) $R$-module has a semilocal endomorphism
ring.
\end{cor}

\begin{proof}
 The assertion follows from Theorem \ref{thetheorem} and Theorem \ref{pr_hdim}.
\end{proof}

\begin{cor}
Let $M$ be a self-projective $R$-module with semilocal endomorphism
ring. Then $\End{M/N}$ is semilocal for any fully invariant
submodule $N$ of $M$.
\end{cor}

\begin{proof}
Since $M$ is self-projective and $N$ fully invariant we get by
\cite[18.2]{wisbauer} that $M/N$ is self-projective. By Theorem
\ref{pr_hdim} we have $hdim(End(M)) = hdim(M) \geq hdim(M/N) = hdim(End(M/N))$.
\end{proof}

Analogous to the fact that a projective module has a semiperfect 
endomorphism ring if and only if it is finitely generated and 
supplemented (see \cite[42.12]{wisbauer}) we get the following corollary:

\begin{cor} \label{cor_lazard}
Let $M$ be a self-projective $R$-module.
$M$ is finitely generated and weakly supplemented if and only if $\End{M}$ is semilocal and $\Rad{M} \ll M$.
\end{cor}

\begin{proof}
 This follows from Theorem \ref{fg_ws}, Theorem \ref{pr_hdim} and the fact
that a module with finite hollow dimension and small radical is finitely
generated.
\end{proof}

The author does not know if the hypothesis of a small radical of $M$ is 
necessary. He raises the following\\

{\bf Question:} Is every (self-)projective 
$R$-module with semilocal endomorphism ring finitely generated ?

\begin{rem} \em 
This question is closely related to an old problem of D.  Lazard. He
considered rings with the property that all projective modules $P$ with
$P/\Rad{P}$ finitely generated are already finitely generated.  Following
 H. Z\"oschinger,  rings with this property are called $L$-rings. He
proved
in \cite{zoeschinger81} that this property is left-right symmetric, i.e.
$R$ is a left $L$-ring if and only if it is a right $L$-ring. Moreover he
showed that a ring $R$ is an $L$-ring if and only if every supplement in
$R$ is a direct summand (\cite[Satz 2.3]{zoeschinger81}). Hence
semiperfect and semiprimitive rings, i.e rings with zero Jacobson radical,
are $L$-rings.  In \cite{jondrup}, S.J\o ndrup showed that every PI-ring is
an $L$-ring. A good resource for some characterizations of $L$-rings
is \cite{mohammed}. 
\end{rem}

\begin{cor} Let $R$ be an $L$-ring and $P$ a projective $R$-module. Then
$P$ is finitely generated and weakly supplemented if and only if $\End{P}$
is semilocal.
\end{cor}

\begin{proof} Assume $\End{P}$ to be semilocal. By Takeuchi's result
(Theorem \ref{pr_hdim}) $P$ has finite hollow dimension and hence
$P/\Rad{P}$ is finitely generated. As $R$ is an $L$-ring, $P$ is finitely
generated.
\end{proof}
\begin{rem} \em
In \cite{gerasimov} V.N. Gerasimov and I.I. Sakhaev constructed 
a non - commutative semilocal ring that is not an $L$-ring (see also
\cite{sakhaev91}, \cite{sakhaev93}). Hence a negative answer to the
question above might be more likely, but the condition of a semilocal
endomorphism ring $\End{P}$ is stronger than $P/\Rad{P}$ being finitely 
generated.
\end{rem}

A ring $R$ is left {\it $f$-semiperfect} or {\it semiregular} if
every finitely generated left ideal has a supplement in $_RR$, equivalently,
 $R/\Jac{R}$ is von Neumann
regular and idempotents in $R/\Jac{R}$ can be lifted to $_RR$ (see
\cite[42.11]{wisbauer}).
Analogous to that we have:

\begin{prop}
For any ring $R$ the following statements are equivalent:
\begin{enumerate}
\item[(a)] every principal left ideal of $R$ has a weak supplement in $_RR$;
\item[(b)] $R/\Jac{R}$ is von Neumann regular;
\item[(c)] every principal right ideal of $R$ has a weak supplement in $R_R$;
\end{enumerate}
\end{prop}

\begin{proof}
$(a) \Rightarrow (b)$ Let $a\in R$. By assumption 
there exists a weak supplement $I \subset R$ of $Ra$.
Hence there exist $b \in R$ and $x\in I$ such that $x=1-ba$. 
Moreover, by Lemma \ref{lem_sem1},
$Rax= Ra\cap Rx \subseteq Ra\cap I \ll R$ implies $ax=a-aba \in \Jac{R}$.
Thus $R/\Jac{R}$ is von Neumann regular.

$(b) \Rightarrow (a)$ For any $a\in R \setminus \Jac{R}$ we get an element
$b\in R\setminus \Jac{R}$ such that $a-aba \in \Jac{R}$. Hence $R(1-ba)$
is a weak supplement of $Ra$ in $_RR$ by Lemma \ref{lem_sem1}. 

$(b) \Leftrightarrow (c)$ analogous.
 \end{proof}

\begin{rem}\em 
It is not difficult to see that $Ra$ has a weak supplement in $_RR$
if and only if $aR$ has a weak supplement in $R_R$ for all $a\in R$. 
The situation for supplements is not that clear. 
Indeed H.Z\"oschinger proved that the property: $Ra$ has a
supplement in $_RR$ implies $aR$ has a supplement in $R_R$ is
equivalent with $R$ being an $L$-ring (see \cite{zoeschinger81}).
\end{rem}

The following proposition, due to S. Page, relates uniform and hollow
dimension.

\begin{prop}[Page] \label{theorem_page}
Let $_RQ$ be an injective cogenerator in $R$-Mod, $T=\End{Q}$ and $M$ any
$R$-module. Then $hdim(_RM) = udim(\Hom{M}{Q}_T)$.
\end{prop}
\begin{proof} See \cite[Proposition 1]{page}. \end{proof}

\begin{rem}\em 
More generally this result can be extended to any injective
cogenerator $_RQ$ in $\sigma[_RM]$ - the full subcategory of $R$-Mod that
contains all $M$-subgenerated left $R$-modules. Hence for any $_RN \in
\sigma[M]$ the formula $hdim(_RN) = udim(\Hom{N}{Q}_T)$ holds where
$T:=\End{Q}$ (see \cite{lomp} for details).
\end{rem}

Using S. Page's result we get another characterization
of semilocal rings in terms of hollow and uniform dimension.

\begin{thm}
The following statements are equivalent for a ring $R$.
 \begin{enumerate}
	\item[(a)] $R$ is semilocal;
	\item[(b)] there exists a generator $_RG$ in $R$-Mod 
		   such that $G_{\End{G}}$ has finite hollow dimension;
	\item[(b')] for any generator $_RG$ in $R$-Mod, $G_{\End{G}}$ has finite hollow dimension; 
	\item[(c)] there exists an injective cogenerator $_RQ$ in $R$-Mod
		   such that $Q_{\End{Q}}$ has finite uniform dimension;
	\item[(c')] for any injective cogenerator $_RQ$ in $R$-Mod, 
                    $Q_{\End{Q}}$ has finite uniform dimension.
\end{enumerate}
In this case $hdim(G_{\End{G}}) = length(R/\Jac{R}) = udim(Q_{\End{Q}})$ 
holds. 
\end{thm}

\begin{proof} 
  Let $_RG \in R$-Mod and $S=\End{G}$.  
  From \cite[18.8]{wisbauer} we know that $_RG$ is a generator in $R$-Mod 
  if and only if $G_S$ is finitely generated, projective in $S$-Mod 
  and $R \simeq End_S(G_S)$. Hence by Theorem \ref{pr_hdim}, we 
  get $hdim(G_S) = hdim(End_S(G_S)) = hdim(R)$. 
  This proves $(a)\Leftrightarrow (b)\Leftrightarrow (b')$. 

  By Page's Proposition \ref{theorem_page} we have $hdim(R) = udim(Q_T)$ 
  where $T=\End{Q}$ for any injective cogenerator $_RQ$.  
  This proves $(a) \Leftrightarrow (c) \Leftrightarrow (c')$. 
\end{proof}

\end{document}